\documentclass{article}
\usepackage[utf8]{inputenc}
\usepackage{appendix}
\usepackage{authblk}
\usepackage{todonotes}

\usepackage{hyperref}

\graphicspath{ {images/} }

\hypersetup{
    colorlinks=true,
    linkcolor=blue,
    filecolor=magenta,      
    urlcolor=cyan,
} 

\title{Parameters for the best convergence of an optimization algorithm On-The-Fly}
\author{Valdimir Pieter $^1$ \\
$^1$ Artificial Intelligence, VU Univeristy Amsterdam}

\date{}

\begin{document}

\maketitle

\section{Plagiarism}
In this paper, the sources is explicitly stated in the text or the footnotes. Moreover, all code used within this paper were designed and implemented from scratch. The goal of this paper is to provide the readers with information that will give them a better understanding on the different types of algorithm and the effect that certain parameters have.

\section{Abstract}
What really sparked my interest was how certain parameters worked better at executing and optimization algorithm convergence even though the objective formula had no significant differences. Thus the research question stated: 'Which parameters provides an upmost optimal convergence solution of an Objective formula using the on-the-fly method?' This research was done in an experimental concept in which five different algorithms were tested with different objective functions to discover which parameter would result well for the best convergence. To find the correct parameter a method called 'on-the-fly' was applied. I run the experiments with five different optimization algorithms. 
\begin{itemize}
    \item Gradient-descent algorithm;
    \item Nelder-Mead algorithm;
     \item Metropolis-Hasting Algorithm;
    \item Simulated Annealing Algorithm;
    \item Evolutionary Algorithm;
\end{itemize}
One of the test runs showed that each parameter has an increasing or decreasing convergence accuracy towards the subjective function depending on which specific optimization algorithm you choose. Each parameter has an increasing or decreasing convergence accuracy toward the subjective function. One of the results in which evolutionary algorithm was applied with only the recombination technique did well at finding the best optimization. As well that some results have an increasing accuracy visualization by combing mutation or several parameters in one test performance. In conclusion, each algorithm has its own set of the parameter that converge differently. Also depending on the target formula that is used. This confirms that the fly method a suitable approach at finding the best parameter. This means manipulations and observe the effects in process to find the right parameter works as long as the learning cost rate decreases over time.


\section{Introduction}
This paper presents how to learn the right parameter for an optimization algorithms. We will have a closer look at five specific ways to optimize a given objective function. This paper will give a closer look at what the affect of these parameters have on the convergence process of optimization algorithm. Eventually what we want to know is what the right parameter is for a given objective function. This means finding the ‘best’ solution among several other solutions.

\begin{itemize}
    \item Gradient-descent algorithm;
    \item Nelder-Mead algorithm;
     \item Metropolis-Hasting Algorithm;
    \item Simulated Annealing Algorithm;
    \item Evolutionary Algorithm;
  
\end{itemize}

 The first algorithm that I will introduce is the Gradient-Descent (GD) algorithm. GD is one of the most popular algorithms to perform optimization. This algorithm uses derivatives of an objective function to find the best solution \cite{ruder2016overview}. GD does this by updating the given input values in the opposite direction of the gradient of the objective function until the minima is found. There are a few different gradient descents.This is called the batch and stochastic gradient descent. They differ by the use of data. For the batch the cost function is the average of the losses. The loss is calculated, one for each data point, based on your prediction. Then, the average of theses losses is taken. However in stochastic gradient descent you calculate your parameter update after each loss. The loss effectively corresponds to the cost. This research will focus solely on the 'Stochastic gradient descent'. Starting with a random point on a function and move in the negative direction of the gradient to reach the local/global minima  \cite{svanberg1981local}.  A research done earlier by Pedregosa states that an appropriate set of hyperparameters is both crucial in terms of model accuracy and computationally challenging. In his work he proposes an algorithm for the optimization of continuous hyperparameters using inexact gradient information. An advantage of this method is that hyperparameters can be updated before model parameters have fully converged \cite{pedregosa2016hyperparameter}.
 He also discussed in his research  that the cost function is costly to evaluate, it is not feasible to perform backtracking line search. To overcome this he used an method in which the step size is corrected according the estimated point from the previous step even though he uses this technique he did not manage to have a formal analysis of the algorithm for this choice of step size \cite{pedregosa2016hyperparameter}.

The second algorithm looked upon is the Nelder-Mead algorithm (NM) also known as the simplex algorithm. It is a gradient-free optimization (GFO) algorithm selected through the method/library in Scipy. It requires function evaluations and is a good choice for simple minimization problems. It uses a triangular shape or simplex to search for an optimal solution \cite{baudin2009nelder}. The simplex shape shifts toward its goal growing shrinking and changing its shape according to a set of rules. Eventually it will convert to the optimal solution \cite{baudin2009nelder}. Nelder-Mead was chosen because it was the most suitable algorithm for mathematical functions and the implementation through scipy was done in (SIAM Journal on optimization) \cite{lagarias1998convergence}.

Metropolis-Hasting (MH), also  known as Monte Carlo method, is used to finding an arbitrary samples sequence from a probability distribution in which direct sampling is difficult to obtain. These sequences are used to approximate the target distribution. MH is a popular way to sample multidimensional distribution, in particular when dimensions are of high value. For lower dimension, there are a few other methods. This paper will however focus on the first. The Metropolis-Hastings algorithm works by generating a sequence of sample values in such a way that, as more and more sample values are produced, the distribution of values more closely approximate the desired distribution. These sample values are produced iteratively, with the distribution of the next sample being dependent only on the current sample value (thus making the sequence of samples into a Markov chain). Specifically, at each iteration, the algorithm picks a candidate for the next sample value based on the current sample value. A research done by Gibson in 1998 stated that: 'An important question for MCMC methods in particular is whether the estimates of this uncertainty obtained from a sequence of samples from the chain is realistic, and whether or not the distributional properties of the generated samples are a true reflection of the target distribution. The results of this paper provide significant evidence that the chains applied do provide a meaningful representation of their equilibrium distributions and that error bounds are realistic' \cite{gibson1998estimating}. My research would like to discover if the parameter will apply to my MH whether it converge to a proper reflection of the target formula.

The fourth sampling algorithm we will look at is the Simulated annealing (SA). SA method is inspired by the process of annealing in metal work. The process goes as follows: the temperature of solid metal is increased until the metal melts. this is done by placing the metal in a heat bath.  when the temperature is high enough the metal melts and becomes liquid. While liquid, the metal is at its the maximum temperature, when high enough its internal structure will change as well as physical properties. Afterwards the metal cools and settle down to reach a very stable flat state. This process is being simulated and is used to generate a solution to distribution sample of the optimization problems \cite{haddock1992simulation}.
A past research done by Pedro A. Castillo stated that a general problem in model selection is to obtain the right parameters that make a model fit the observed data. In this study it was crucial to find the appropriate weights and learning parameters when dealing with models that have a Multilayer perception that is trained with back propagation \cite{castillo1999sa}. Their strategy  attempts to avoid Lamarckism  but, at the same time, it is a good strategy
to avoid local minima. 

Lastly, Evolutionary algorithm(EA) \cite{fonseca1995overview} , is a nature inspired approach to optimization. Which is the process of getting the most out of something and making it better during this process. We are in the search for the best or optimal solution to a problem. EA is inspired by the idea of survival of the fittest from Darwinian evolution and Mendel's modern genetics. Such as reproduction, mutation, recombination, and selection. Candidate solutions to the optimization problem play the role of individuals in a population, and the fitness function determines the quality of the solutions \cite{lan2019evolving}. Evolution of the population then takes place after the repeated application of the above operators. The general idea behind the algorithm is that if biological evolution can produce something as amazing as humans over many generations then we should be able to use the same process artificially to evolve optimal solutions for vehicles \cite{xu2019online}, aircraft \cite{lan2016UAV}, spaceships and robots \cite{lan2018directed,lan2019evolutionary,lan2019simulated}. To dispel any over hyping it must be emphasized that in practice evolutionary computation offers approximations of optimal solutions to difficult problems. Often it can be over hyped as the all-powerful algorithm, where every problem in existence is a nail  however this is not the case. EA consist of the following steps:
\begin{itemize}
    \item (Init) Initialize a population of solutions and evaluate;
    \item (Generate) Generate new solutions by applying mutation and  recombination;
    \item (Evaluate) Evaluate new solutions;
    \item (Select) Replace the least-fit individuals of the population with new individuals;
  
\end{itemize}
A similar test was performed for EA in a paper written by Eiben in 2017 in which he adjusted the setting parameters on-the-fly \cite{eiben2007parameter}.  This research will have the same approach for all five of the algorithms. It is also stated in the same paper that 'The other motivation for controlling parameters on-the-fly is the assumption that the given parameter can have a different “optimal value in different phases of the search. If this holds, then there is simply no
optimal static parameter value; for good EA performance one must vary this
parameter' \cite{eiben2007parameter}. In the final section of this research I hope to validate this statement based on my own results. Therefore, the research question states: What parameters provide an upmost best convergence solution of a given objective function when applying the 'on-the-fly' method?

All these different methods used by other research above proves that none of them are perfect and works only on certain algorithm. With my analysis, I try to prove that the on-the-fly method should perform well enough to apply to any objective function converge with the right found parameters without the need to consider which Algorithm is applied. Another contribution is that these finding can be further experimented with to add value to the already existing research in the many different ways to find the right set of parameters.


\section{Problem statement}
In the following section the problem set for each algorithms and the usage of the objective function that correspond with each problem will be discussed. Specify what the objective function is and comment on possible difficulties that are expected to encounter. There will be five different target fomula/objective functions with each algorithm having a set of repetitions. The results of these test runs will be discussed in a later section.

\subsection{Gradient descent and Nelder-Mead}
As mentioned in the introduction, optimization is choosing inputs that will result in the best possible outputs. This can mean a variety of things from deciding on the most effective allocation of available resources, producing a design with the best characteristics to choosing control variables that will cause a system to behave as desired.

Optimization problems often involve words such as maximize or minimize. Optimization is also useful when there are limits or constraints on the resources involved or boundaries
restricting the possible solutions \cite{lan2018ICARCV}. For our problem, we went with an objective function to where we had to find the minima. With domain restriction for x1 and x2:

\begin{equation}
    \int(x) = x_1^2+2x_2^2-0.3cos(3\pi x_1)-0.4cos(4\pi x_2)+0.7 
\end{equation}
\begin{equation}
    (x_1,x_2) \in [-100,100]^2
\end{equation}

The objective function for Nelder-Mead:
\begin{equation}
    \int(x) = (x_1 + 2x_2 - 7)^2 + ( 2x_1 + x_2 -5)^2 
\end{equation}
\begin{equation}
    (x_1,x_2) \in [-100,100]^2
\end{equation}
Minima, is another word for the lowest point of a function. The goal is to find this minima with use of Gradient descent. While finding this minimum we tried diffent values, inputs and compare the algorithms against each other if needed. For some, this would have been a simple problem, because at times it is easy to see the correct solution, but for more complicated problems it can be difficult to immediately see the correct solution. More detail of the optimization formula can be found in the paper by Chen published in 2017  \cite{chen2017well,lan2020time}.

Guessing and checking may be time consuming. This can result in an expansive cost to run an algorithm. Furthermore, it can be difficult to find the correct parameter such as alpha, step-size and number of iterations which also have effect on the quickness to find the minima. So using the on-the-fly mechanic \cite{eiben2007parameter} will be a perfect way to approach this problem. Another difficulty is to implement the algorithm due to many ways to program a code. That being said, just finding the basic solution for the met minimum is a difficulty on its own. some functions can have many local minima which make it harder to find the lowest minimum coordination \cite{svanberg1981local}.

\subsection{Metropolis-Hasting and Simulated Annealing}
Optimization is also useful when there are limits or constraints on the resources involved or boundaries, restricting the possible solutions. This is what makes optimization useful. For our problem, we were given a Target Distribution to where we had to sample by using MH and SA algorithm. The objective function that is used for Metropolis-Hasting is:

\begin{equation}
\small
p(x)\propto exp\{-0.01(sin(x_1)exp[(1-cos(x_2))^2]+cos(x_2)exp[(1-sin(x_1))^2]+(x_1-x_2)^2)\}
\end{equation}

\begin{equation}
    (x_1,x_2) \in [-3,3] \times [2,4]
\end{equation}

The objective function for Simulating Annealing is: 
\begin{equation}
\small
p(x)\propto exp\{-0.02(cos(x_1)exp[(1-sin(x_2))^2]-sin(x_2)exp[(1+cos(x_1))^2]-(x_1-x_2)^2)\}
\end{equation}

\begin{equation}
    (x_1,x_2) \in [-3,3] \times [2,4]
\end{equation}

through sampling, we can replace the analytical calculation. Instead, we generate two candidate points using a proposal distribution. Continue with evaluating the candidate points and eventually either reject or accept the new  generated candidate points based on the acceptance probability.

\begin{figure}[!ht]
    \centering
    \includegraphics[width=0.6\textwidth]{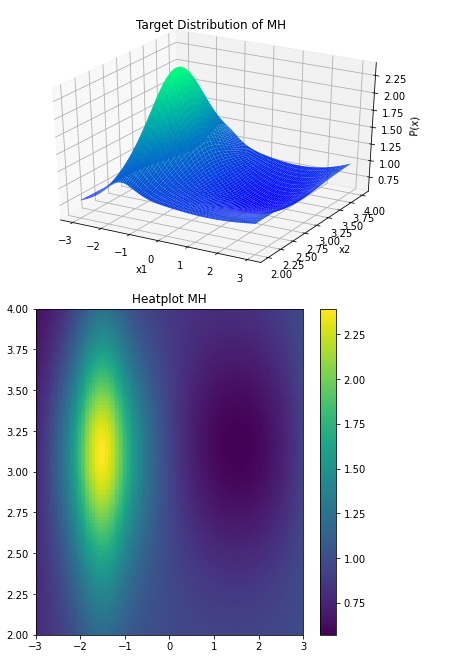}
    \caption{Objective Function Graph MH}
    \label{fig:mhsapic1}
\end{figure}

\begin{figure}[!ht]
    \centering
    \includegraphics[width=0.6\textwidth]{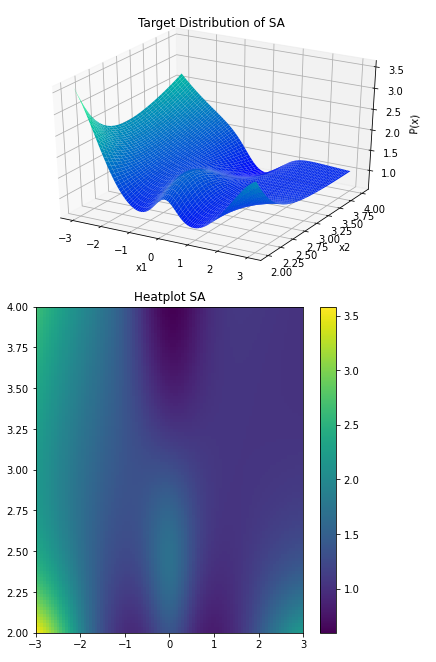}
    \caption{Objective function of SA}
    \label{fig:mhsapic11}
\end{figure}

One of the problems that occur while performing a sampling algorithm is the amount of points being sampled. A Large number of sampling points will take much longer to execute and may cause oversampling \cite{cho2018sampling}. This can be very expansive and realistically not reachable. However, having to little points will lead to under-sampling and cause a very inaccurate distribution of the target distribution. By adjusting different values, hyper parameters and iterations I expect to find the right values for the algorithm to perform efficiently and accurate.

\subsection{Evolutionary Problem Statement}
In the introduction I stated that the fitness would determine how well the candidate scores points. The goal is to implement components of an evolutionary algorithm: a recombination operator, a mutation operator and selection mechanisms, and analyze their behavior. I am interested in optimizing a given the function Repressilator \cite{tomczak2020differential}. that could be queried, but the gradient input cannot be calculated. The repressilator is a genetic regulatory network consisting of at least one feedback loop with at least three genes, each expressing a protein that represses the next gene in the loop. In biological research, repressilators have been used to build cellular models and understand cell function. The input to the system is a vector of: 

\begin{equation}
    x=[\alpha_0,\eta,\beta,\alpha]^\tau \in [-2,10] \times [0,10] \times [-5,20] \times[500,2500]
\end{equation}

The biggest limitation of EA is that it cannot guarantee optimally. The solution quality also deteriorates with the increase of problem size. However it can generate good quality solutions for any problem and function type \cite{lan2020learning}. Stochastic algorithms in general can have difficulty obeying equality constraints. the EA is sensitive to the initial population used. Wide diversity of feasible solutions is optimal in such case. In the paper written by Tomczak in \cite{tomczak2020differential} he concluded in the following quotation that: 'We provide a theoretical analysis of the proposed linear operators by proving their reversibility, and inspecting their eigenvalues. Further, we show empirically on three testbeds (benchmark function optimization, discovering parameter values of the gene repressilator.  \cite{tomczak2020differential}
systems, and learning neural networks) that producing new candidates on-the-fly allows to obtain better results in fewer number of
evaluations compared to DE'
This shows that the On-The-Fly mechanic seems to work well for different type of algorithms.


\section{Methodology}
\subsection{GD and NM}
\paragraph{}To have our data-set to be useful and understandable I performed different types of methods. Gradient descent is a (machine) learning algorithm that is used to draw a fit a line to a set of points to find the best minimum solution.
and it is useful initially, I calculated the partial derivatives analytically. This was needed to eventually find the next best step towards the minimum.
below you find the two partial derivatives that were calculated analytic:

\begin{equation}
    f\;'(x_1)=2x_1+0.9\pi*sin(3\pi x_1)
\end{equation}

\begin{equation}
    f\;'(x_2)=4x_2+1.6\pi*sin(4\pi x_2)
\end{equation}

\begin{figure}[!htb]
    \centering
    \includegraphics[scale=0.50]{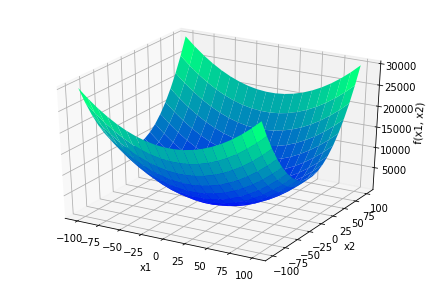}
    \caption{Objective Function Graph of GD}
    \label{fig:gddfpic1}
\end{figure}

\begin{figure}[!htb]
    \centering
    \includegraphics[scale=0.75]{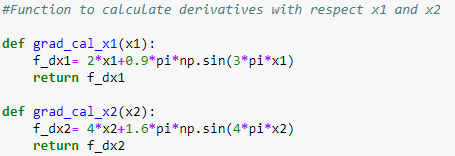}
    \caption{Derivatives code}
    \label{fig:gddfpic2}
\end{figure}

\begin{figure}[!htb]
    \centering
    \includegraphics[scale=0.45]{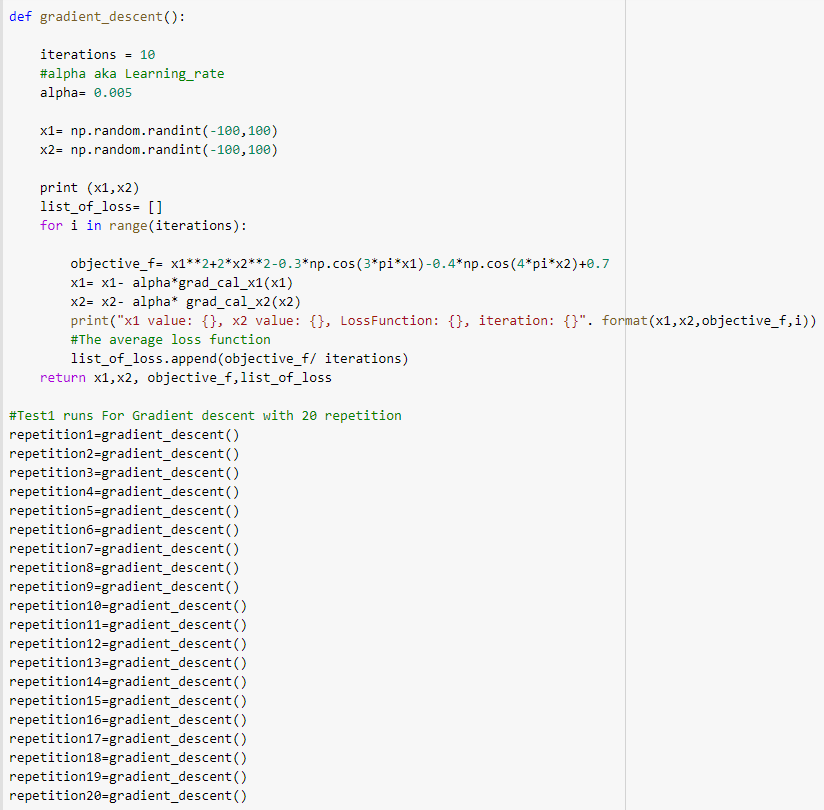}
    \caption{Gradient descent code and Repetitions}
    \label{fig:gddfpic3}
\end{figure}

\begin{figure}[!htb]
    \centering
    \includegraphics[scale=0.5]{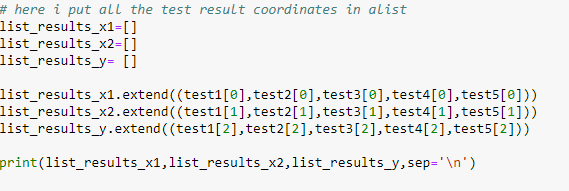}
    \caption{List of returned Values}
    \label{fig:gddfpic4}
\end{figure}

The following step was to implement the objective function and GD in python code. Starting with the import of packages; Math, Numpy, Matplotlib and Scipy. These packages helped with translate the function in python language. I coded the objective function to eventually visualize and have a better understanding of how the  graph looked like. Which is showed in Figure~\ref{fig:gddfpic1}. For the two derivatives I made two separate 'DEF functions' that returned the partial derivatives as seen in Figure~\ref{fig:gddfpic2}. In addition, I made a third 'def-function' that calculated the Gradient-Descent(Figure~\ref{fig:gddfpic3}).
within the function I placed the Iteration, Alpha and the random input values for x1 and x2.
After I got the results from the first test run consisting of 20 repetitions I inserted all the results in a new list divided by x1, x2 and y values. This would make the visualization less complex. These results are visualised through scatter boxes, tables and graphs.(Figure ~\ref{fig:gddfpic4})

\begin{figure}[!ht]
    \centering
    \includegraphics[scale=0.50]{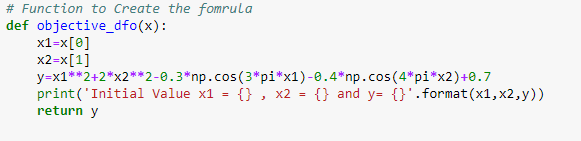}
    \includegraphics[scale=0.40]{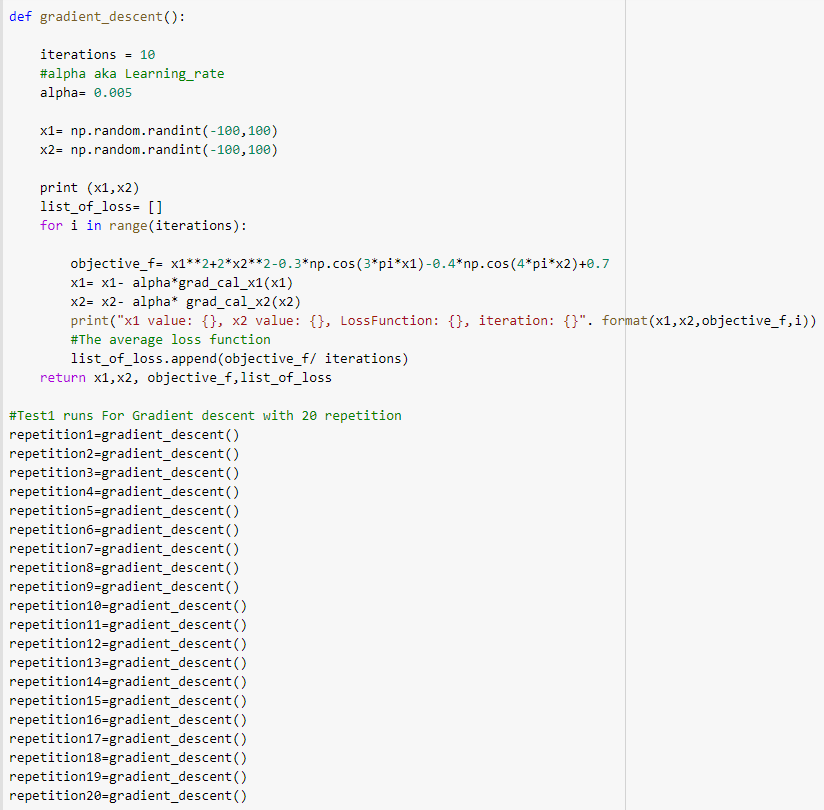}
    \caption{Nelder-Mead Code}
    \label{fig:gddfpic5}
\end{figure}

For the derivative-free algorithm(DFO), Nelder-Mead. I used Scipy package to get the optimization implemented. It was not to difficult to code this due to a large portion of the algorithm calculation was pre-coded. All I had to do was feed it an objective function and (hyper)parameters and let the algorithm run Test 2. As done with GD algorithm, with NM I also did several repetitions, each time changing one value whilst rest maintained consistent and observe  correlation or informative chance. In the next chapter, we will have a detailed overview of the manipulated parameter and values provided~\ref{fig:gddfpic5}

\subsection{MH and SA }
\begin{figure}[!ht]
    \centering
    \includegraphics[scale=0.80]{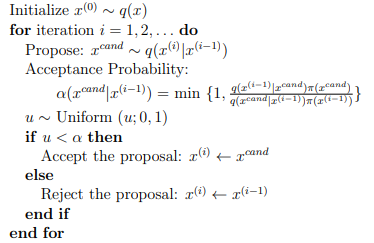}
    \caption{Pseudo code of If statement}
    \label{fig:mhsapic2}
\end{figure}

For Metropolis-Hasting, I had the mathematical notations translated in python language to continue with the implementation of the target distribution. Figure~\ref{fig:mhsapic1} shows the plot of the target distribution for Metropolis-Hasting and Figure~\ref{fig:mhsapic11} the plot for Simulating Annealing along with the corresponding heat-plots. With the use of Heat-plots, we can define the volume of points and mark the important areas efficiently.
The MH algorithm is slightly divided into three parts: Generate, evaluate and selecting. For starters, an arbitrary starting value was picked and then iterative accepting or rejecting candidate samples drawn from another distribution. For obtaining a proposal distribution point I did so by drawing random samples from a normal (Gaussian) distribution. A Random distribution was chosen because a proposal distribution is a symmetric distribution. According to Yildirim(2012) choices of symmetric proposals include Gaussian distributions or Uniform distributions centred at the current state of the chain  \cite{yildirim2012bayesian}. This means that when the proposal distribution is symmetric the acceptance probability becomes proportional to how likely each of the current state and the proposed state are under the full joint density.
The second step was evaluating the candidate point generated by proposal distribution. This was done in a separate function in which the Alpha was calculated as the acceptance probability. After calculating the alpha it was checked whether the candidate point should be accepted or rejected through an IF-statement conditions. If alpha was greater than the standard deviation the candidate point would be accepted and if it was smaller than the candidate point it would be rejected. The pseudo-code can be seen in Figure~\ref{fig:mhsapic2}. 

\begin{figure}[!ht]
    \centering
    \includegraphics[scale=0.45]{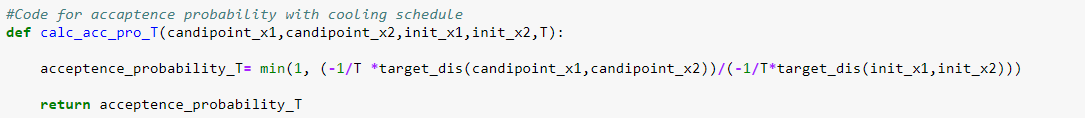}
    \caption{Temperature influenced Acceptance Probability}
    \label{fig:mhsapic3}
\end{figure}

For the Simulated Annealing algorithm  a major part of the code was reused for the implementation. However SA was centered around sample through cooling schedule. What was an addition to the code was the Temperature influence. The alpha was recalculated with influence of the Temperature(T). With each iteration the temperature is updated as seen in figure~\ref{fig:mhsapic3}.

\subsection{Evolutionary}

Evolutionary algorithm consist of a few stages \cite{fonseca1995overview}, the first stage of any evolutionary algorithm is the initialization. We can think about this stage as some kind of Big Bang where we create our initial population out of random generator(unifrom).
We call a group of solutions to a problem at any stage of the progress a population. It is important to define our data structure often referred to as a chromosome in the context of evolutionary algorithms. This chromosome simply stores a set of items or genes which describe our solutions. We can think of this as some type of blueprint.
With the initial population generated we move on to the evaluation stage. This is where we figure out the fitness of each of our solutions. The Fitness is often a single number which gauges how good a solution is. The higher the Fitness number the better the solution. We need this information later on to apply this idea of survival of the fittest. Of course this Fitness value depends entirely on the problem we are trying to solve. Do do this we need the right parameters set. Therefore is the use of the 'on-the fly' method our way to tackle this. 

The next stage is called selection and this is where survival of the fittest
takes place. If we were to pick three solutions from this population to produce offspring for future generations. Which three would we pick? typically we would select the Three solutions with the highest Fitness score because we expect these three to contain the most valuable genetic information and we want this genetic information to be passed on to the new solutions.
Finally, variation stage  we will aim to use our selected solutions which we can refer to as the parent population to produce new solutions in which we can refer to as the offspring population. One method of variation is called chromosome recombination here you can see the same chromosome. The crossover points mix two parent solutions to produce a single offspring solution. Another method of variation which is typically used together with crossover is called mutation. With mutation all we do is select a random gene on the chromosome
of a randomly selected solution then we add a random number to this gene. This introduces new genetic information into our population. This is important as it allows us to explore different types of solution which may perform better on our problem.

\section{Experiments}
In this following section I will discuss how the experiments on each algorithm is performed in detail. Note that parameter are picked on-the-fly initially and as the algorithms proceed these will change accordingly.
\subsection{GD and NM experiment}
Firstly we will look at the values given for the GD algorithm:
Both x1 as x2 were given random starter value between the constrain range from inclusive -100 to 100. These were set to be the hyper parameter for alpha= 0.001 and iteration of 10.
A smaller value was chosen to start with to observe the type of results or influence these parameters had on the loss cost. Twenty repetitions were done to find the lowest loss cost. Only the initial start point would differ with each repetition. Based on the results of the loss I would slightly increase or decreases the Iteration first. While observing to see if the loss function would lower over time. Was this the case then the I would proceed with adjusting the Alpha. Would The Alpha size lead to a drastically change of the loss function I would adjust the Alpha smaller value compare to the previous value it had. 
Important is to know that 'Alpha' means the same as step-size or learning-rate. I wanted discover the right Alpha and iteration parameter as mention above by increasing and decreasing until loss function reached the smallest possible outcome. Thus the convergence of the GD would be successfully executed. This same 'on-the-fly' method I would apply to all five different test scenarios. 

For Nelder-Mead algorithm, x1 and x2 had a random initial value generated like GD.  Alpha is called here ‘atol’. For NM this is 0.005 and Maxiter is the maximum allowed number of iterations and function evaluations. This was set on 100 to start with.
With this experiment, I expected to find the local minima as well to find the best alpha and iterations while the remaining values was kept constant. Only the input would be different each run.  A randomly generated input each run would give a different starter point on the function. Moreover, based on the actual iteration NM needed to find the smallest loss function Alpha would be re-evaluated and adjusted if needed.

\subsection{MH and SA experiment}
For the metropolis-hasting I started with the  following parameters:
\begin{itemize}
    \item Initial value for x1 and x2 were  2 , 2 respectively.;
    \item Iteration (N) of  1000;
    \item Standard deviation of 0.2;
  
\end{itemize}

 The initial values were generated randomly and placed in the proposal distribution function, which generated the candidate points.  However, the initial points as the generated candidate points were checked if they met the constrains domain given. If this was not the case a  new candidate point would be generated. This was done with a specific function. This function used the Gaussian distribution, it is a normal distributed with a standard deviation of 0.2.
 
 Afterward we calculated the acceptance probability with the following formula:
 \begin{equation}
    \frac{(1, (-target\_dis(candipoint\_x_1,candipoint\_x_2))}{(-target\_dis(init\_x_1,init\_x_2)))}
\end{equation}

After generating the acceptance probability/alpha which consisted of the candidate points and initial points. Now we had to check whether these candidate points should be accepted or not. To do so I had a hyper parameter 'u' which was  randomly generated  value through a uniform distribution between 0 and 1. With each iteration this ‘u’ value would change. Comparable to the acceptance probability that would be calculated with each iteration. This is put in a the table of Figure~\ref{fig:mhsapic5}.

Furthermore, with each value I appended these candidate points to a rejected or acceptance list based on the outcome. Also kept count of the cumulative number of accepted points which starts at 0. If a candidate point is (of was) accepted, a point would be added. If the old point was picked, one point would be subtracted. This experiment was done over 20 repetitions with a larger iteration values and completely new initial value for x1 and x2 which was set to -3  and 2 respectively. The other values remained the same.  The goal for MH was to see what kind of effect the amount of iteration had on the results. I expect to find a convergence that will look like the the target distribution with each time the iteration size gets bigger.

For the Simulated Annealing algorithm experiment I used a slightly adjusted version of the code from MH. This was done for a better comparison between to algorithms. I also started with the iteration value of N=100 However, what I did adjust was the acceptance probability formula calculation. The initial temperature value (T) and and the cooling(C) schedule was applied. The cooling schedule which was a constant value that decreased the Temperature over iterations.(Figure~\ref{fig:mhsapic3})
Parameters for Simulated annealing:

\begin{itemize}
    \item Initial value for x1 and x2 were  -3 , 2 respectively.
    \item Iteration (N) = 1000  and  10000
    \item C = 1
    \item T= 100
    \item $T= C*0.95^i$
\end{itemize}

The Simulated annealing experiment was ran for 20 repetitions, with iteration of 10 and C, a constant of 1.
The initial value for x1 and x2 was -3  and 3  respectively. The other values remained the same.
The goal for SA was to discover whether there was a significant change between the repetitions when changing parameters such as iteration to a much higher value. Also adjusting the Temperature to manipulate its effect on  the sampling points.

\subsection{Evolutionary experiment}
Lastly, For EA experiment I went with set run of repetitions. What is important to know is that the first few repetitions had the same pattern. Only the fourth test run differentiated. By default, the first few repetitions had:  
\begin{itemize}
    \item Population size = 100
    \item Generations size = 10
    \item Standard deviation = 1

\end{itemize}
these values remained constant through the three test runs. This was critical for the a better comparison of results when applying solely the mutation, permutation or both functions. For detailed explanation these methods view the paper written by Spears in \cite{spears1998role}. 

At first, the initialization occurred and gave a randomly generated x values that were between the constraints given. Continue giving these value to the EA class the fitness was calculated. These x values were also given to recombination, mutation or both depending on the operators given. Repetition one  x, f =EA.step(x,f,0,1) indicates that running in class EA calling the function step with operators x values fitness of x, no recombination (0) and mutation (1). For the second test run, I did the same but instead of mutation I  only applied recombination. Which looked as follows : x,f= EA.step(x,f,1,0)

For the third run again the same population size and generation, the standard deviation was used but this time we run both mutation and recombination at the same time. 
x,f EA.step(x,f,1,1)  indicating  1 for applying recombination and second 1 for mutation. The results of these tests will be discussed in the next chapter

For the fourth and last run, I wanted to experiment with different population size. so I picked a population of 15 and a standard deviation of 10  and called the step function with both the mutation and recombination.
\begin{itemize}
    \item Population size = 15
    \item Generations size = 10
    \item Standard deviation = 0.5
\end{itemize}

\section{Results and Conclusion}

\subsection{Gradient descent and Nelder-Mead results}
\begin{figure}[!ht]
    \centering
    \includegraphics[scale=0.6]{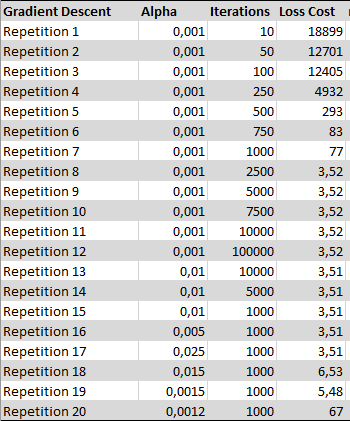}
    \caption{Gradient Descent result table}
    \label{fig:gddfpic6}
\end{figure}

Figure ~\ref{fig:gddfpic6} shows the table of the Gradient Descent. Each  Alpha, Iteration and lowest Loss cost found per repetition. The first twelve repetitions have an Alpha  of 0.001. Only the iteration value were adjusted over these first 12 repetitions. Furthermore, in repetition 15 the Loss cost was on its lowest. the table also shows that over time the Alpha was adjusted through repetition 12 through 17 which resulted in the Loss costs to be equal. The goal was to find the best alpha value and iterations.  
As the Loss cost values are decreasing over each iteration. 

\begin{figure}[!ht]
\begin{minipage}{0.5\linewidth}
\includegraphics[width=\textwidth]{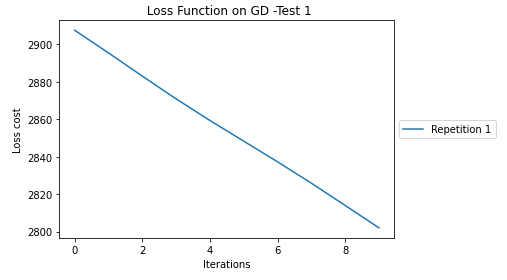}
\caption{Repetition 1 Loss Cost curve }
\label{fig:gddfpic7}
\end{minipage}%
\hfill
\begin{minipage}{0.5\linewidth}
\includegraphics[width=\textwidth]{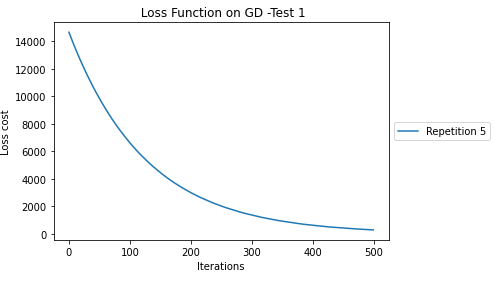}
\caption{Repetition 5 Loss cost curve}
\label{fig:gddfpic8}
\end{minipage}%
\hfill
\begin{minipage}{0.5\linewidth}
\includegraphics[width=\textwidth]{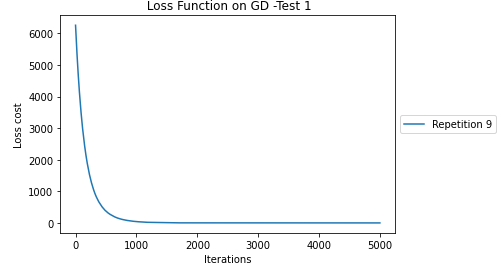}
\caption{Repetition 9 Loss cost curve}
\label{fig:gddfpic9}
\end{minipage}
\begin{minipage}{0.5\linewidth}
\includegraphics[width=\textwidth]{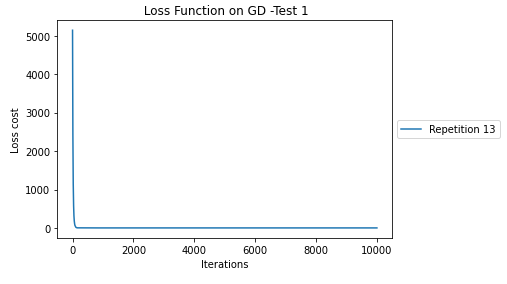}
\caption{Repetition 13 Loss Cost curve }
\label{fig:gddfpic10}
\end{minipage}%
\hfill
\begin{minipage}{0.5\linewidth}
\includegraphics[width=\textwidth]{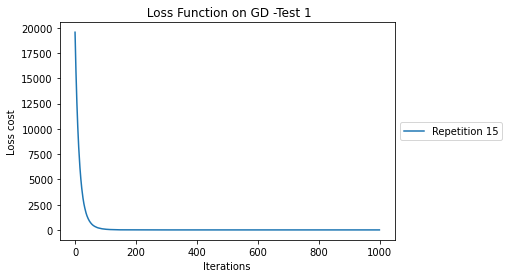}
\caption{Repetition 15 Loss cost curve}
\label{fig:gddfpic11}
\end{minipage}%
\hfill
\begin{minipage}{0.5\linewidth}
\includegraphics[width=\textwidth]{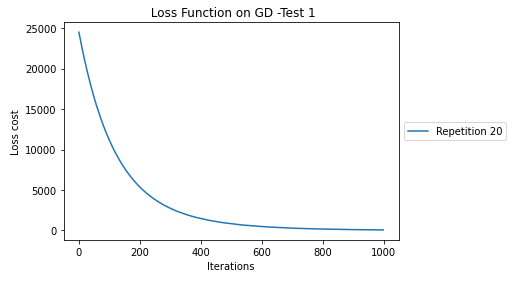}
\caption{Repetition 20 Loss cost curve}
\label{fig:gddfpic12}
\end{minipage}
\end{figure}

In Figures ~\ref{fig:gddfpic7} to ~\ref{fig:gddfpic12} the process of the adjustment of the iteration and alpha shows the effect the parameters have on the Loss cost. You can also see that the decreasing of the Loss cost is not consistent. Which means that I had to adjust iteration one more time to see if it could get the Loss cost to decrease more toward zero. Repetition 9 with an iteration of 5000 has a much deeper Loss cost curve than repetition 5. Based on this result the decision was made to stop increasing the iteration and place the focus to find the right Alpha parameter as  seen in Figure ~\ref{fig:gddfpic12}. This figure shows that the Loss cost becomes slightly smaller than the first 13 repetitions. 
 
\begin{figure}[!ht]
    \centering
    \includegraphics[scale=0.60]{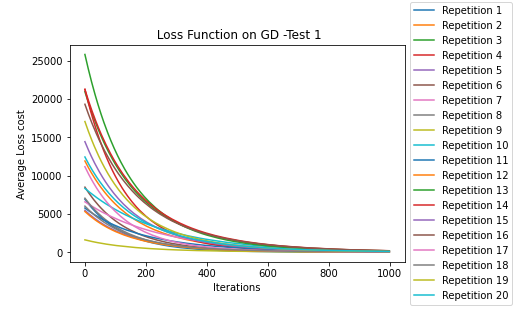}
    \caption{All 20 repetition Loss cost }
    \label{fig:gddfpic13}
\end{figure}

 The smaller the iteration size the less expansive it is to execute an algorithm  as researched by Buyya in \cite{buyya2005scheduling}  and a  to small Alpha may lead to over sample. Therefore  an Alpha 0.001 and iteration of 1000 is chosen as the best possible solution in Test scenario one. With nuanced conclusion
 In Figure ~\ref{fig:gddfpic13} the average loss cost of 20 repetitions is visualized with corresponding best Alpha and iteration chosen 0.001 and 1000. 
 
\begin{figure}[!ht]
\centering
\includegraphics[scale=0.6]{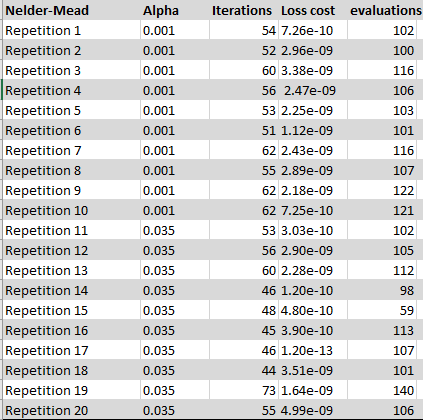}
\caption{Nelder-Mead Test 2 Result table }
\label{fig:gddfpic14}
\end{figure}

The results for Nelder-Mead algorithm, Test scenario 2 is shown in table of Figure~\ref{fig:gddfpic14}. The best solution found was repetition 14 which had the lowest loss cost. However, the amount of iterations that occurred for to find this Loss cost was not the best possible outcome, instead repetition 18 found  the (local) Loss cost in iteration of just 44. Even though the best solution had a much higher iteration value.  Further more the the scatter plot shows that four out of the five runs had a very low minima.

When comparing the two Algorithm against each other you can see that NM algorithm performs much better than GD algorithm.  The size of iterations needed to find the best solution is significant lower than the GD algorithm. while GD needed 1000 iteration to find the lowest Loss cost NM only needed 46. Also the overall Loss cost values of NM are generally lower than the Lost cost from GD. 

Lastly the standard deviation of  0.19and mean of 0.09 of Nelder-Mead  smaller compared to Gradient Descent mean of 1.86 and standard deviation of 1.05 which means that there is little to no variation for Nelder-Mead. Therefore, Nelder-Mead is a much better optimization algorithm than Gradient Descent based on these results given.

\begin{figure}[!ht]
        \centering
        \includegraphics[width=0.6\textwidth]{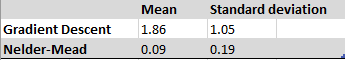}
        \caption{The mean and Standard deviation table }
        \label{fig:gddfpic15}
\end{figure}

\subsection{Metropolis-Hasting and Simulated Annealing results}

\begin{figure}[!ht]
    \centering
    \includegraphics[scale=0.5]{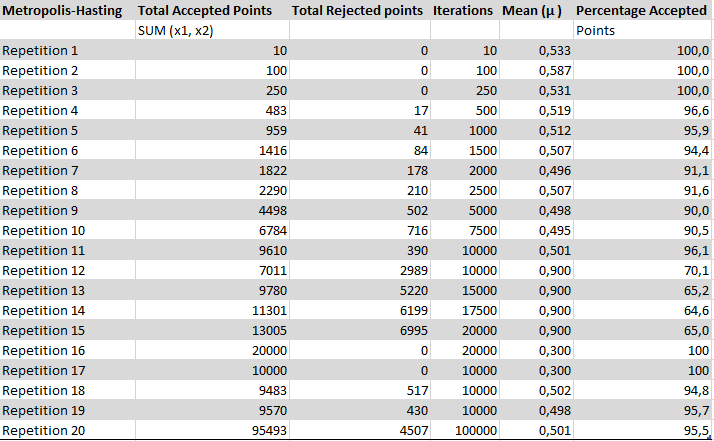}
    \caption{Result Table of Metropolis-hasting }
    \label{fig:mhsapic5}
\end{figure}

\begin{figure}[!ht]
    \begin{minipage}{0.5\linewidth}
    \includegraphics[width=\textwidth]{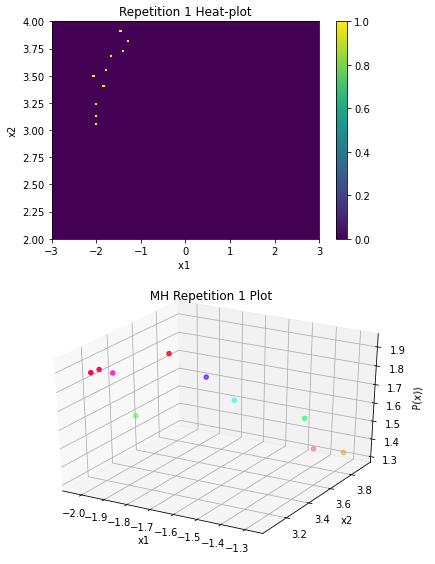}
    \caption{Repetition 1: 3D graph MH}
    \label{fig:mhsapic6}
    \end{minipage}%
    \hfill
    \begin{minipage}{0.5\linewidth}
    \includegraphics[width=\textwidth]{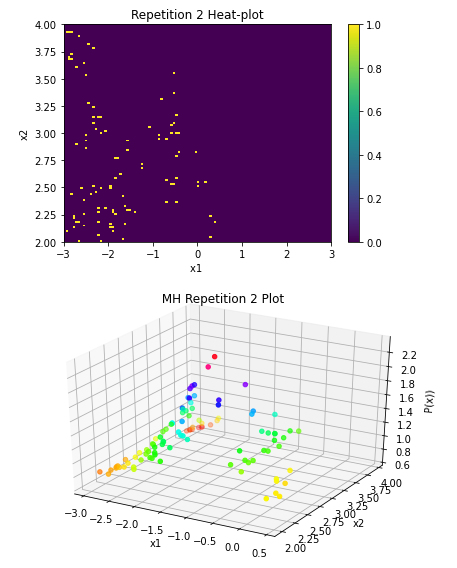}
    \caption{Repetition 2: 3D graph MH}
    \label{fig:mhsapic7}
    \end{minipage}%
\end{figure}

\begin{figure}[!ht]
\begin{minipage}{0.5\linewidth}
\includegraphics[width=\textwidth]{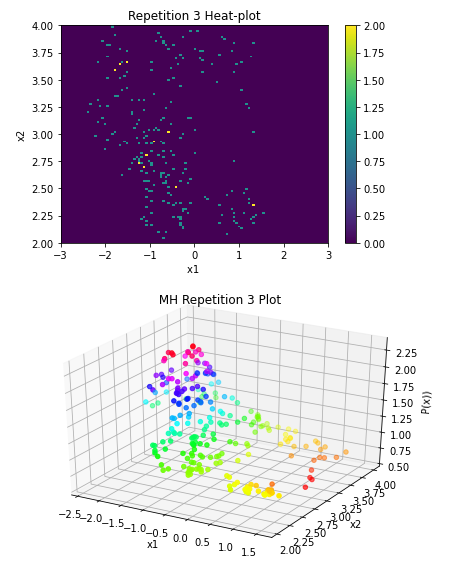}
\caption{Repetition 3  3D graph MH}
\label{fig:mhsapic8}
\end{minipage}%
\hfill
\begin{minipage}{0.5\linewidth}
\includegraphics[width=\textwidth]{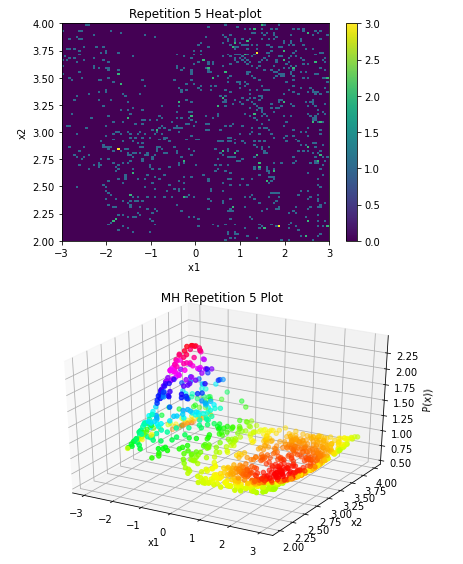}
\caption{Repetition 5: 3D graph MH}
\label{fig:mhsapic9}
\end{minipage}%
\end{figure}

\begin{figure}[!ht]
\begin{minipage}{0.5\linewidth}
\includegraphics[width=\textwidth]{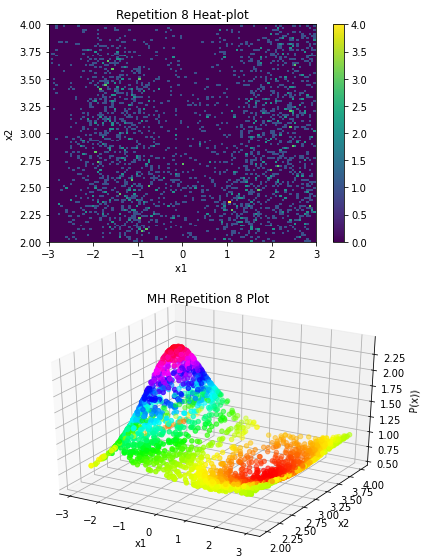}
\caption{Repetition 8  3D graph MH}
\label{fig:mhsapic10}
\end{minipage}%
\hfill
\begin{minipage}{0.5\linewidth}
\includegraphics[width=\textwidth]{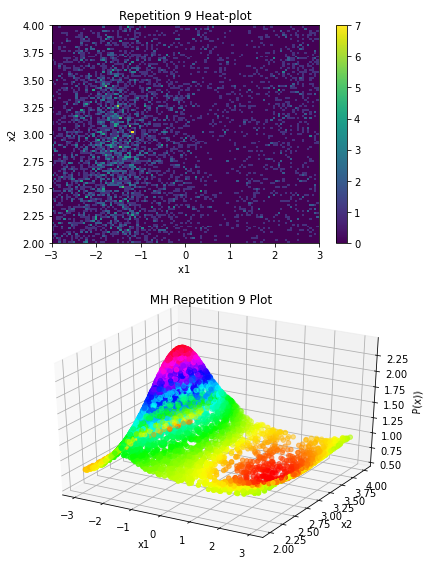}
\caption{Repetition 9  3D graph MH}
\label{fig:mhsapic110}
\end{minipage}%
\end{figure}

\begin{figure}[!ht]
\begin{minipage}{0.5\linewidth}
\includegraphics[width=\textwidth]{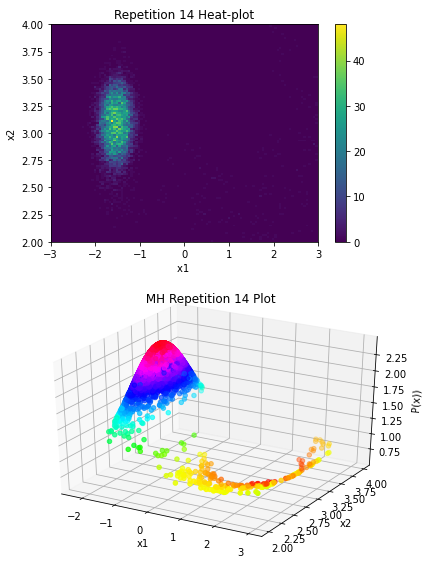}
\caption{Repetition 14  3D graph MH}
\label{fig:mhsapic12}
\end{minipage}%
\hfill
\begin{minipage}{0.5\linewidth}
\includegraphics[width=\textwidth]{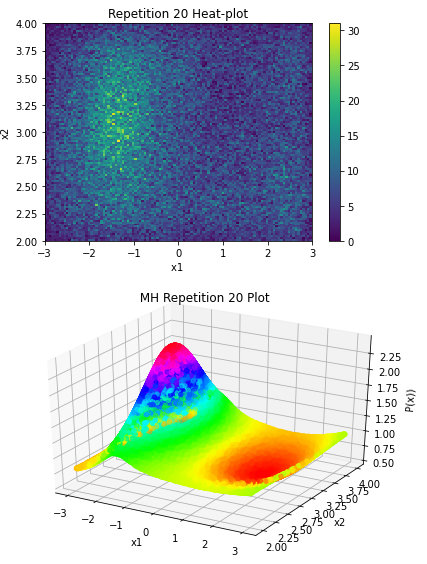}
\caption{Repetition 20  3D graph MH}
\label{fig:mhsapic13}
\end{minipage}%
\end{figure}

Figure~\ref{fig:mhsapic5} shows the table of Test Three. The table presents the results of the total accepted and rejected points as well as the iteration size, percentage of the accepted point and the mean of (u)- acceptance probabilities. Repetition 1 through 3 shows an accepted percentage of 100. Although the iterations are rather small and does not provide a proper convergence of the objective function as shown in figure~\ref{fig:mhsapic6}. 
Furthermore, the table shows that repetition 13 and 14 have the lowest accepted percentage. The cause of this relies on the acceptance probability that was chosen. A lower acceptance probability accept more points this can be seen in the repetitions in which the acceptance probability is 0.900 (column 4). Knowing this provides a honest deliberation on choosing the right iteration parameter now that is proven that acceptance probability around 0.5 preforms well enough in Test Three. Figure~\ref{fig:mhsapic6}  through ~\ref{fig:mhsapic13} shows the 3D graphical representation of the convergence of the objective function. As the the size of the iteration becomes larger the similarity of the object function becomes noticeable. Also the heat-plots above each graph plot represents tensity of the correlation between the candidates points x1 and x2 of the accepted points. The lighter the color, the higher the accepted candidate points tensity. As seen in Figure ~\ref{fig:mhsapic12} repetition 14 has dominantly more points towards 40 tensity than ~\ref{fig:mhsapic10}.

\begin{figure}[!ht]
    \centering
        \includegraphics[scale=0.5]{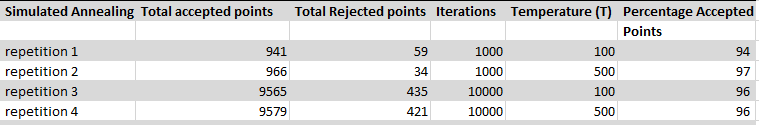}
    \caption{The First 4 Repetition of SA}
    \label{fig:mhsapic14}
\end{figure}

\begin{figure}[!ht]
\begin{minipage}{0.5\linewidth}
\includegraphics[width=\textwidth]{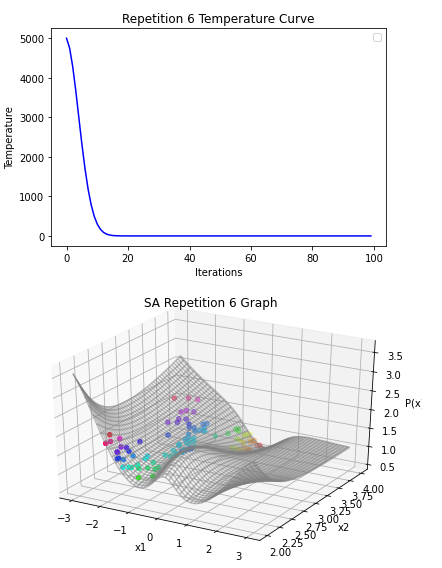}
\caption{Repetition 6  3D graph SA}
\label{fig:mhsapic15}
\end{minipage}%
\hfill
\begin{minipage}{0.5\linewidth}
\includegraphics[width=\textwidth]{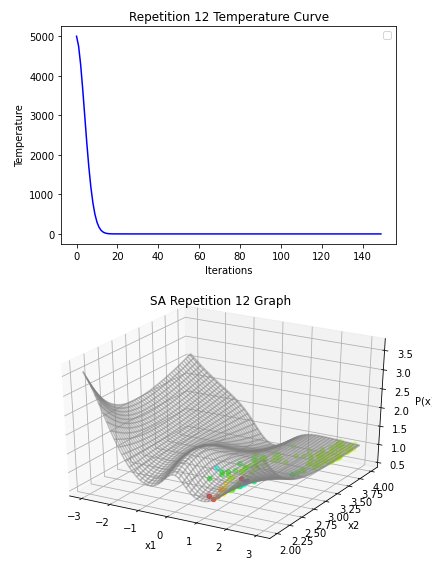}
\caption{Repetition 12  3D graph SA}
\label{fig:mhsapic16}
\end{minipage}%
\end{figure}

\begin{figure}[!ht]
\begin{minipage}{0.5\linewidth}
\includegraphics[width=\textwidth]{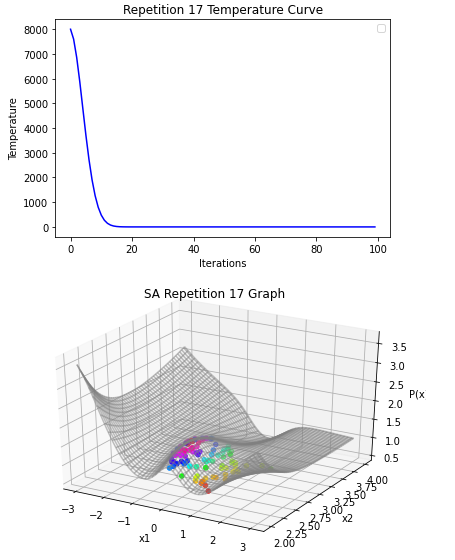}
\caption{Repetition 17  3D graph SA}
\label{fig:mhsapic17}
\end{minipage}%
\hfill
\begin{minipage}{0.5\linewidth}
\includegraphics[width=\textwidth]{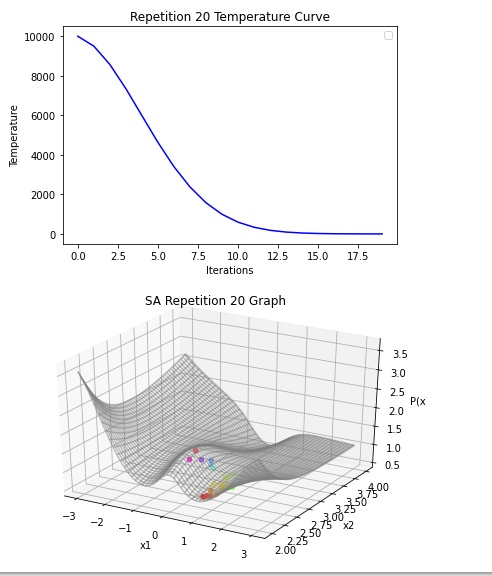}
\caption{Repetition 20  3D graph SA}
\label{fig:mhsapic18}
\end{minipage}%
\end{figure}

Simulated Annealing Test run 4 results gives that the the temperature applied for the first four repetition does not have a significant effect on the accepted candidate points (Figure ~\ref{fig:mhsapic14}). more over,  If we first look at the repetition 1 with an iteration size of 1000 we can see that when we increase the temperature from 100 to 500 while constant remains the same the amount of accepted points are better at a temperature of 500. The results shows us again that the accepted points are slightly better than T of 100. Also a higher Temperature sample does not necessarily converge to a much accurate result. If we look at the Temperature curve from repetition 6,12,17 and 20(Figure ~\ref{fig:mhsapic15} to ~\ref{fig:mhsapic18}. It shows that around 15 iteration the temperature does not decrease as swiftly. Therefore the iteration is adjusted to a lower size in repetition 20. 
Here you can see that majority of the accepted points fall on the left side of the heat plot. Just like in the target distribution Figure ~\ref{fig:mhsapic1}. The more iteration happens the more point are located on the left side. Which is converge to accurate representation of the target distribution. 

As we would have expected the two traditional algorithms repeatedly stick at points associated with singularities. looking at the solution corresponding to the graphical , we find that both SA and MH algorithms generally converged to a similarity of their objective function. Both shows high acceptance percentage with parameters for iterations  1000 and alpha/acceptance probability of 0.5.

\subsection{Evolutionary results}

\begin{figure}[!ht]
    \centering
    \includegraphics[scale=0.5]{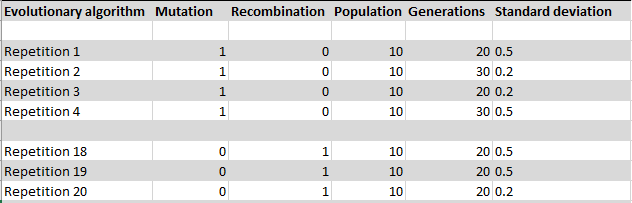}
    \caption{Evolutionary Algorithm Repetition table}
    \label{fig:evopic1}
\end{figure}

\begin{figure}[!ht]
\begin{minipage}{0.5\linewidth}
\includegraphics[width=\textwidth]{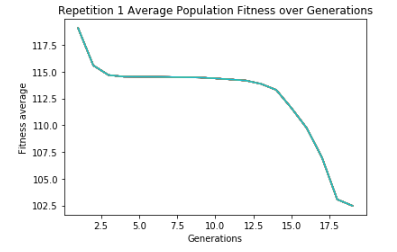}
\caption{EA-Repetition 1}
\label{fig:evopic2}
\end{minipage}%
\hfill
\begin{minipage}{0.5\linewidth}
\includegraphics[width=\textwidth]{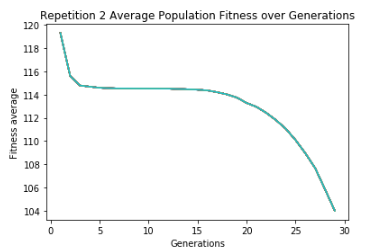}
\caption{EA-Repetition 2}
\label{fig:evopic3}
\end{minipage}%
\hfill
\begin{minipage}{0.5\linewidth}
\includegraphics[width=\textwidth]{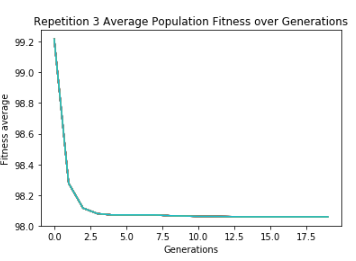}
\caption{EA-Repetition 3 }
\label{fig:evopic4}
\end{minipage}%
\hfill
\begin{minipage}{0.5\linewidth}
\includegraphics[width=\textwidth]{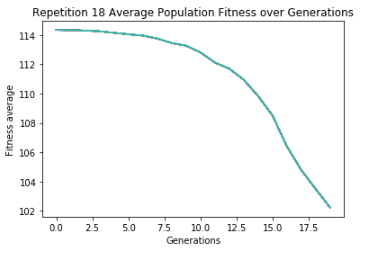}
\caption{EA-Repetition 18 }
\label{fig:evopic5}
\end{minipage}%
\hfill
\begin{minipage}{0.5\linewidth}
\includegraphics[width=\textwidth]{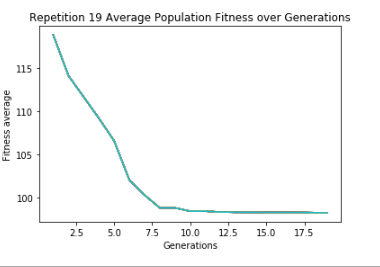}
\caption{EA-Repetition 19 }
\label{fig:evopic6}
\end{minipage}%
\begin{minipage}{0.5\linewidth}
\includegraphics[width=\textwidth]{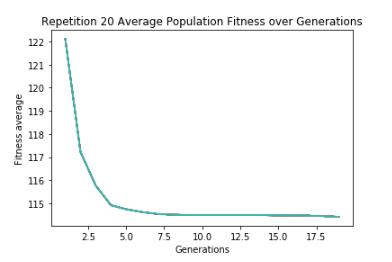}
\caption{EA-Repetition 20 }
\label{fig:evopic7}
\end{minipage}%
\end{figure}

Figure ~\ref{fig:evopic1} shows the results of Test 5 as set in a table. Repetitions 1-10 had mutation applied whilst repetition 11-20 had recombination. It shows that the population size remained constant as the focus was on the generation and standard deviation(std). The average mean fitness reached its lowest point in the 19th repetition (Figure~\ref{fig:evopic2}). Recombination was applied here with std of  0.5. Furthermore, the average mean fitness converge with quickness between Generation 1 to 5 in which the curves decreases the steepest among  the 20 repetitions. However over all the repetitions done, non of the curves got close enough to the zero to claim whether the chosen parameters performed well enough. What it does show is that the std used in the experiment can be set to a larger size. Based on the results showing a very small to no significant convergence to 0 between the repetitions. This could be a good exploration idea for the next research. In addition, the methods in this paper can be used to the optimization in aircraft \cite{lan2016developmentuav}, knowledge graph \cite{liu2020influence}, sensor network \cite{lan2015bayesian,lan2016development,lan2017development}.

Overall, with the results discussed above, we can conclude that parameters chosen on-the-fly are a good way to find the right parameter but not the 'best'. it does take, some repetitions to actually find some correlation. At the same time each optimization algorithm parameter has its effect on the result of the convergence process. Whether the parameter has significant or rather small differences. For the Gradient descent the iteration over 1000 converged to a perfect similarly of the objective function while the parameters of evolutionary algorithms, would depend on the specific problem of recombination and mutation. Multiple repetitions of the algorithms with different probability e.g. 0.5, 00,1. The parameters also depend on the implementation of the code. One would tend to use smaller Alpha/Step-size rates for Simulating annealing compared to deviate free algorithms for example Nelder-Mead.

\bibliographystyle{unsrt}
\bibliography{main}

\end{document}